# Super-Linear Sub-Quadratic Worst-Case Computational Complexity Achieved by *spdspds* Algorithm for Linear Programming Problem

© Dr.(Prof.) Keshava Prasad Halemane,
Professor - retired (2017JAN31) from
Department of Mathematical And Computational Sciences,
National Institute of Technology Karnataka Surathkal,
Srinivasnagar, Mangaluru - 575025, India.
https://orcid.org/0000-0003-3483-3521

Independent Researcher
(no affiliations; no funding; no data used; no conflict of interests)
SASHESHA, 8-129/12 Sowjanya Road, Naigara Hills,
Bikarnakatte, Kulshekar Post,
Mangaluru-575005. Karnataka State, India
https://www.linkedin.com/in/keshavaprasadahalemane/
k.prasad.h@gmail.com

## Abstract

The Symmetric Primal-Dual *Sym*plex Pivot Decision Strategy (*spdspds*) is a novel iterative algorithm to solve linear programming problems. A symplex pivoting operation is considered simply as an *exchange* between a basic (dependent) variable and a non-basic (independent) variable, in the *Goldman-Tucker Compact Symmetric Tableau* (CST) which is a unique symmetric representation common to both the primal as well as the dual of a linear programming problem in its standard canonical form. From this viewpoint, the *classical simplex pivoting* operation of Dantzig may be considered as a restricted special case.

The *infeasibility index* associated with a symplex tableau is defined as the sum of the number of primal variables and the number of dual variables that are infeasible. A *global effectiveness measure* of a pivot selection is defined/determined as/by the *decrease in the infeasibility index* associated with such a pivoting operation. The selection of the symplex pivot element is made by seeking the best possible anticipated *decrease in the infeasibility index* from among a wide range of candidate choices with non-zero values. The worst-case computational complexity of the *spdspds* algorithm is shown to be $O(L^{1.5})$ where L refers to the problem-size expressed in terms of the size(length) of the input data.

## 1. Introduction

Linear Programming (LP) problem represents one of the most widely used class of computational models, for which any possible improved solution technique would certainly be highly desirable. Suggested alternative solution strategies include variations in the *classical simplex method* of Dantzig [1], later followed by the Ellipsoid Method of Khachiyan [2] [3] and the Karmarkar Algorithm [4] - both classified now as belonging to Interior Point Algorithms. Terlaky [5] Todd [6] and Adler et al [7] present an overview of the various developments. Chinneck [8] gives an overview of the developments in the use of the concept of *feasibility and infeasibility in optimization problems.* Let us not get diverted much into the historical developments etc.

The simplex pivoting operation of Dantzig represents a move from one basis/tableau to another basis/tableau, by/through a *single exchange* between an entering *infeasible* non-basic variable and a leaving *feasible* basic variable. For a chosen entering *infeasible* non-basic variable, the leaving *feasible* basic variable is to be selected so as to meet certain restrictive criteria in terms of the corresponding limitation of moving only between two *neighboring extreme/vertex points* of the polytope defined by the set of linear system of inequalities - so as to maintain feasibility while improving the objective function value by moving further towards the optimum.

The proposed *spdspds* approach can be considered as a novel generalization of the simplex method of Dantzig, in terms of lifting all of such restrictions and providing a wider scope for the selection of the pivots - any nonzero element of the coefficient matrix in the tableau can be a potential candidate pivot element. It is indeed true that the very term *simplex pivot* has been redefined here - as a simple/single *exchange* between an entering non-basic variable and a leaving basic variable - maintaining only the combinatorial/structural property of being a simple/single *exchange* between a selected pair - that being the justification for renaming it as *symplex pivot* - emphasizing the primal-dual symmetry therein. The symplex pivot of *spdspds* need not necessarily correspond to a pair consisting of an entering *infeasible* non-basic variable and a leaving *feasible* basic variable; also it does not require to be limited to a move between neighboring extreme/vertex points of the associated polytope - although the move does indeed correspond to one between a pair of *intersection points* defined by the set of linear system of inequalities.

The actual selection of a *spdspds* pivot element is governed by an analysis of the associated measure of goodness of such a pivot choice. A global measure of goodness or a *global* *effectiveness measure* (*gem*) for pivot selection is defined, utilizing the novel concept of *infeasibility index* associated with a symplex tableau - defined as the sum of the number of primal variables and the number of dual variables that are infeasible. The change in the infeasibility index associated with a symplex pivot element can be determined by a thorough analysis of the tableau data. To guarantee the best computational performance, it is proposed to select a pivot element corresponding to the best possible decrease in the infeasibility index. Note that the concept of *infeasibility measure* [9] was introduced, *for the first time*, in the context of *algorithm development* for *optimization under uncertainty* [10] by this author in his doctoral research work at CMU during 1979-1982. This paper is derived from the latest version of full-length preprint [21] of this research work; that has been made available in several public open access repositories.

It will be shown that the spdspds algorithm passes through a non-repeating sequence of CST tableaus and reaches a terminal tableau wherefrom no further reduction in the infeasibility index is possible. The length of such sequence is limited by the infeasibility index of the initial tableau, except for cases of degeneracies that may possibly cause certain elongated sequence. An analysis of the data pattern in the *terminal tableau* can be used to classify the problem into one of the possible *six categories* - an *infeasibility index* of *zero* indicates optimum solution. Even in the absence of an optimum solution, the *spdspds* algorithm allows one to explore/determine the most suitable alternative solutions, including a comprehensive parametric analysis, etc.

## 2. Goldman-Tucker Compact-*Sym*metric-Tableau CST

We will go through some well-known preliminaries for the sake of establishing the notational conventions used in this report, as used in our earlier reports [11] [12] [13] [14] on spdspds.

The Symmetric Primal-Dual Pair of LP in the Standard Canonical Form [15] is as follows:

Primal Problem:

$$\begin{aligned} \text{maximize} \quad & c.x = f \\ \text{s.t.} \quad & A.x \leq b \\ & x \geq 0 \end{aligned} \tag{1}$$

Dual Problem:

$$\begin{aligned} \text{minimize} \quad & v.b = g \\ \text{s.t.} \quad & v.A \geq c \\ & v \geq 0 \end{aligned} \tag{2}$$

The descriptions for each of the problem parameters in (1) & (2) above are as follows:

| | | |
|---|---|---|
| x | Primal decision variables | n x 1 vector |
| c | Primal objective function coefficients | 1 x n vector |
| f | Primal objective function value | 1 x 1 scalar |
| A | Primal constraint coefficient matrix | m x n matrix |
| b | Primal constraint upper bound | m x 1 vector |
| v | Dual decision variables | 1 x m vector |
| g | Dual objective function value | 1 x 1 scalar |

We introduce the m x l vector y of slack variables to (1) and the 1 x n vector u of surplus variables to (2) to rewrite the symmetric primal-dual pair in an equivalent formulation as follows:

Primal Problem:

$$\begin{array}{llcl} \text{maximize} & c.x + 0.y & = & f \\ \text{s.t.} & A.x + I_m.y & = & b \\ & x,\ y & \geq & 0 \end{array} \qquad (3)$$

Dual Problem:

$$\begin{array}{llcl} \text{minimize} & v.b + u.0 & = & g \\ \text{s.t.} & v.A - u.I_n & = & c \\ & v,\ u & \geq & 0 \end{array} \qquad (4)$$

This Symmetric Primal-Dual pair is represented in the Goldman-Tucker Compact *Sym*metric Tableau (CST) as shown in Figure-1.

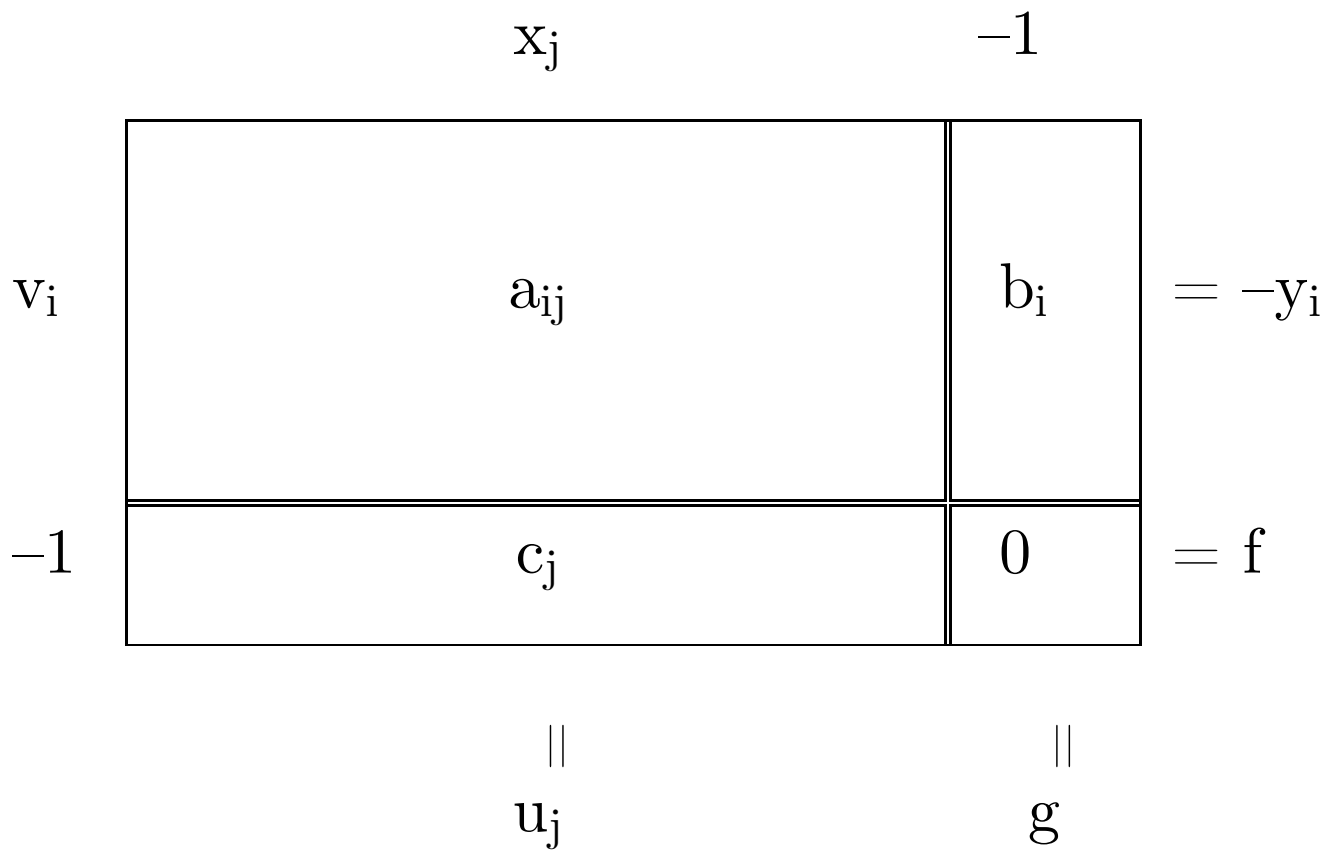

Figure-1: Initial Compact *Sym*metric Tableau ($T_0$)

The relations in (3) & (4) corresponding to the Primal and Dual can be combined to give -

$$c.x + u.x \quad = \quad v.A.x \quad = \quad v.b - v.y \qquad (5)$$

$$f \;\; + u.x \quad = \quad v.A.x \quad = \quad g \;\; - \;\; v.y \qquad (6)$$

For any feasible (basic or non-basic) solution to the P-D pair,

$$u,\ x,\ v,\ y \qquad \geq \qquad 0 \qquad (7)$$

and therefore,

$$(g - f) \quad = \quad (u.x + v.y) \quad \geq \quad 0 \qquad (8)$$

For any basic solution (feasible or infeasible) to the P-D pair,

the non-basic variables are set to zero;

that is,

$$u,\ y \qquad = \qquad 0 \qquad (9)$$

and therefore,

$$(g - f) \quad = \quad (u.x + v.y) \quad = \quad 0 \qquad (10)$$

The entries in the Goldman-Tucker Compact-*Sym*metric-Tableau (CST) directly correspond to the associated basic solution to the P-D pair, thus establishing a one-to-one correspondence between a solution basis and the associated Goldman-Tucker Compact-*Sym*metric-Tableau (CST). A pivoting operation on the CST tableau corresponds to the associated move from one basis to another.

For the LP problem P-D pair (1) & (2) or equivalently (3) & (4) the above tableau in Figure-1 represents the initial tableau indicating the *initial basic solution* (IBS) wherein $y_i$ are the primal basic variables associated (one-to-one permanent association) with $v_i$ the dual non-basic variables, whereas $x_j$ are the primal non-basic variables associated (one-to-one permanent association) with $u_j$ the dual basic variables. Note that the column-labels $x_j$ (along with the label –1 for the RHS column) and the row-labels $v_i$ (along with the label –1 for the objective function row) play a significant role in the CST tableau, and the way to interpret (read) the CST tableau is as follows:

Primal Problem:

$$\sum_{j \in C} a_{ij} \,.\, x_j - b_i = -y_i \,, \quad i \in R \ \text{(row index)}$$

$$\sum_{j \in C} c_j \,.\, x_j - 0 = f \qquad \text{(function to be maximized)} \tag{11}$$

Dual Problem:

$$\sum_{i \in R} v_i \,.\, a_{ij} - c_j = u_j \,, \quad j \in C \ \text{(Column index)}$$

$$\sum_{i \in R} v_i \,.\, b_i - 0 = g \qquad \text{(function to be minimized)} \tag{12}$$

wherein the variables $x_j$, $y_i$, $v_i$, $u_j$ are all considered to be non-negative.

## 3. Algebra (Arithmetic) of the *Sym*plex Pivoting Process

With the Goldman-Tucker Compact-*Sym*metric-Tableau (CST) representation for linear programming, in its standard/canonical form, one can observe that once a pivot element is selected, the actual pivoting process (the *algebra* and hence the *arithmetic* operations) is the same irrespective of the pivot selection; for example whether it is a primal pivot or a dual pivot. Hence it suffices to present here a single (common) set of operations representing the actual pivoting process - be it primal or dual. This *expressional elegance* and *computational efficiency* along with the *convenience* and the *versatility* (as will be evident later) are the reasons why the above representation has been selected for the purpose of our study.

For the sake of generality, let us imagine that we are somewhere in the middle of solving a LP problem (say after the $k^{th}$ iteration), and have the system model represented by a tableau ($T_k$) as shown in Figure-2.

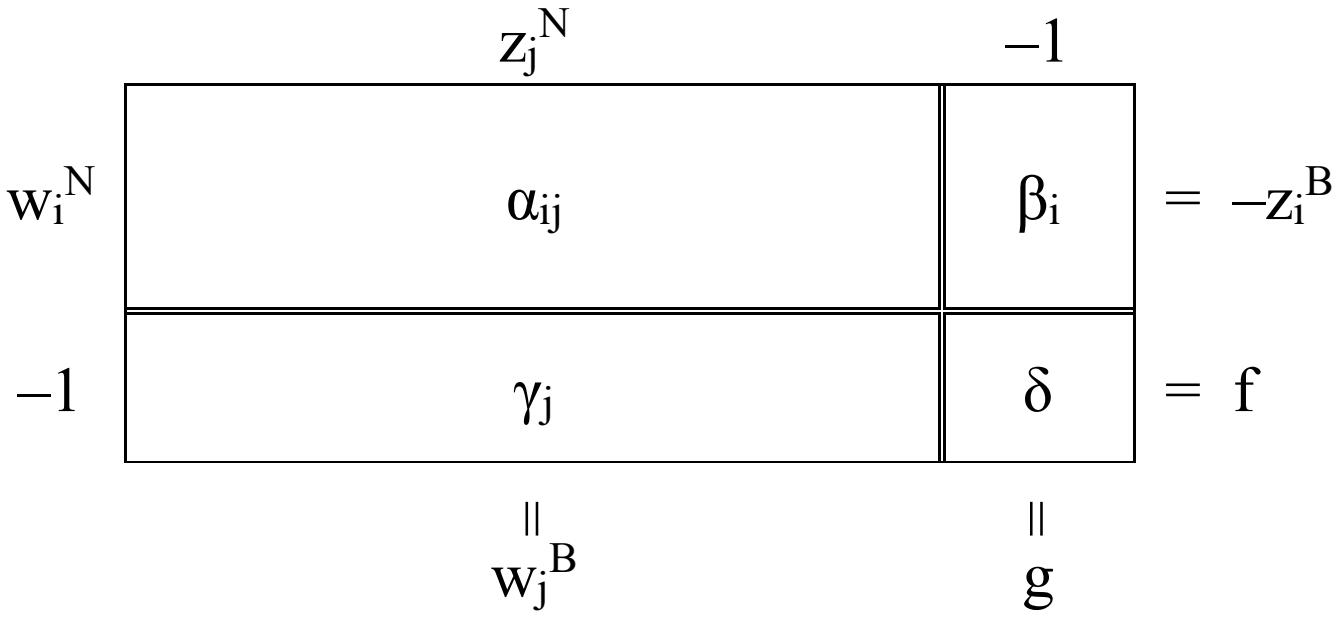

Figure-2: Compact *Sym*metric Tableau ($T_k$) after k iterations

By the nature of the sequence of elementary row (column) operations being performed during any pivoting process, the system model represented in Figure-2 is equivalent to that represented by the initial tableau which corresponds to the P-D pair (11) & (12). The transformed version of the primal-dual pair directly expressed by the above tableau is as follows:

Primal Problem:

$$z_i^B = \beta_i - \sum_{j \in C} \alpha_{ij} \cdot z_j^N , \qquad i \in R \ \text{(row index)}$$

$$f = -\delta + \sum_{j \in C} \gamma_j \cdot z_j^N , \qquad \text{(function to be maximized)} \tag{13}$$

Dual Problem:

$$w_j^B = -\gamma_j + \sum_{i \in R} w_i^N \cdot \alpha_{ij} , \qquad j \in C \ \text{(column index)}$$

$$g = -\delta + \sum_{i \in R} w_i^N \cdot \beta_i \qquad \text{(function to be minimized)} \tag{14}$$

The effect of a pivoting operation on (13) & (14) performed with a chosen pivot element $\alpha_{IJ}$ is exactly to affect an exchange between the P-D variable pairs indicated by I and J in (13) and (14). That is, $z_J^N$ is entered into primal basis in exchange for $z_I^B$ in (13), and $w_I^N$ is entered into dual basis in exchange for $w_J^B$ in (14). Suppose we have chosen the pivot element $\alpha_{IJ}$ using some appropriate pivot selection scheme, and we would like to derive the resulting tableau ($T_{k+1}$). Let the resulting tableau ($T_{k+1}$) be indicated in Figure-3.

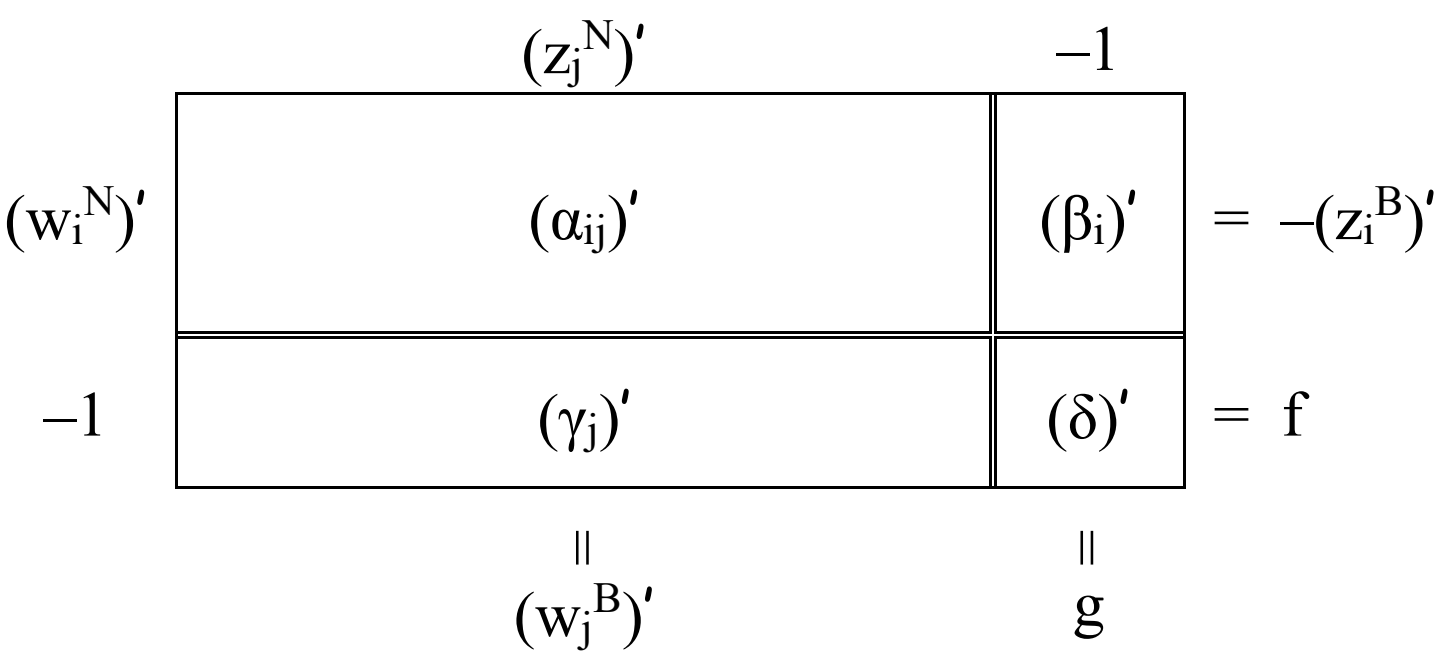

Figure-3: Compact *Sym*metric Tableau ($T_{k+1}$) after (k+1) iterations

The algebra (*arithmetic*) of deriving the above tableau ($T_{k+1}$) of Figure-3 from the previous tableau ($T_k$) of Figure-2 is detailed below (equally applicable to the generic spdspds "*omni*plex" pivoting operation mentioned later):

| | | |
|---|---|---|
| $(\alpha_{IJ})' \leftarrow (\alpha_{IJ})^{-1}$ | $(\alpha_{Ij})' \leftarrow (\alpha_{IJ})^{-1}\alpha_{Ij}$ | $(\beta_I)' \leftarrow (\alpha_{IJ})^{-1}\beta_I$ |
| $(\alpha_{iJ})' \leftarrow -\alpha_{iJ}(\alpha_{IJ})^{-1}$ | $(\alpha_{ij})' \leftarrow \alpha_{ij} - \alpha_{iJ}(\alpha_{IJ})^{-1}\alpha_{Ij}$ | $(\beta_i)' \leftarrow \beta_i - \alpha_{iJ}(\alpha_{IJ})^{-1}\beta_I$ |
| $(\gamma_J)' \leftarrow -\gamma_J (\alpha_{IJ})^{-1}$ | $(\gamma_j)' \leftarrow \gamma_j - \gamma_J (\alpha_{IJ})^{-1}\alpha_{Ij}$ | $(\delta)' \leftarrow \delta - \gamma_J (\alpha_{IJ})^{-1}\beta_I$ |

*along with an exchange of label*s associated with pivot row I and pivot column J; that is effectively:

$$(z_J^N)' \leftarrow z_I^B; \quad (z_I^B)' \leftarrow z_J^N; \quad (w_J^B)' \leftarrow w_I^N; \quad (w_I^N)' \leftarrow w_J^B;$$

while retaining the very same labels for all the other rows and columns; that is:

$$(z_j^N)' \leftarrow z_j^N; \quad (z_i^B)' \leftarrow z_i^B; \quad (w_j^B)' \leftarrow w_j^B; \quad (w_i^N)' \leftarrow w_i^N;$$

for $i \in R\backslash\{I\}$ and $j \in C\backslash\{J\}$.

It is to be noted here that the Goldman-Tucker Compact-*Sym*metric-Tableau (CST) is a *unique symmetric* representation *common* to both the primal as well as the dual of a linear programming problem in its standard canonical form. Also, the tableau evolves from $T_0$ as the initial tableau representing the problem in its standard canonical form, while following the sequence of symplex iterations all the way to the final terminal tableau $T_*$, while always being a tableau representation of an *equivalent system* of linear inequalities along with the corresponding objective function.

From (13) & (14) above, we can get:

$$(f+\delta) + \sum w_j^B \cdot z_j^N = \sum \sum w_i^N \cdot \alpha_{ij} \cdot z_j^N = (g+\delta) - \sum w_i^N \cdot z_i^B \qquad (15)$$

and therefore,

$$(g - f) \;=\; \sum w_i^N \,.\, z_i^B \;+\; \sum w_j^B \,.\, z_j^N \qquad (16)$$

The values of the primal and the dual basic variables as well as (the primal & the dual) objective function value corresponding to a basis/tableau can be directly read from the entries of the tableau - the primal objective function value being always the same as the dual objective function value for every basis/tableau all along the sequence of symplex iterations - each of the summation terms in both (15) and (16) above being zero, and the *complementary slackness condition being automatically satisfied for every basis/tableau.* The relationships (15) & (16) above indicate the effect of moving away from a basis, as the non-basic variables are moved away from zero; in particular, one can note from (16) above that *for all feasible solutions, the dual objective function value is an upper bound for the primal objective function value* and *the primal objective function value is a lower bound for the dual objective function value*, which is a well known relationship. A basic solution for the P-D pair is optimal if & only if it is basic feasible for the P-D pair.

It will be clear from the foregoing discussion that any specific pivoting operation will affect/change the feasibility of a primal-basic-variable [dual-basic-variable] if the ratio $\beta_i/\alpha_{iJ}$ is closer to zero than the *pivoting-ratio* $\beta_I/\alpha_{IJ}$ [if the ratio $\gamma_j/\alpha_{Ij}$ is closer to zero than the *pivoting-ratio* $\gamma_J/\alpha_{IJ}$]: the feasibility/infeasibility of those rows [columns] with the associated ratios farther away from zero beyond (relative to) the corresponding *pivoting-ratio* will not get affected by such pivoting operation.

## 4. Infeasibility Index : A Global Effectiveness Measure

As an *inverse measure of goodness*, the infeasibility index $\lambda$ of a given CST tableau is defined as the sum of the *primal infeasibility index* $\mu$ and the *dual infeasibility index* $\nu$. It corresponds to the number of basic variables in primal & dual which are infeasible in the given tableau. That is, we define $\lambda = \mu + \nu$; and therefore -

$$\lambda = (\mu, \text{number of rows with } \beta < 0) + (\nu, \text{number of columns with } \gamma > 0);$$

or

$$\lambda = (\mu, \text{number of rows with } z_i^B < 0) + (\nu, \text{number of columns with } w_j^B < 0).$$

*If the infeasibility index* $\lambda$ *of the given tableau equals to* zero *then it indicates that the tableau is the terminal tableau which is feasible and optimal. By the definition of the infeasibility index, it* can never be negative*, nor can it be more than the sum of the number of columns and the number of rows in the Compact Symmetric Tableau.* That is,

$$0 \leq \lambda = \{(0 \leq \mu \leq m) + (0 \leq \nu \leq n)\} \leq (m + n).$$

Given a tableau, the anticipated *change in the infeasibility index* ($\tau = \Delta\lambda$) can be associated with each cell that can be a potential candidate pivot element (that is, $\alpha \neq 0$). This change in the infeasibility index consists of two components, one is the *change in the primal infeasibility index* ($\sigma = \Delta\mu$) and the other is the *change in the dual infeasibility index* ($\rho = \Delta\nu$). That is,

$$\text{Change in the infeasibility index, } \tau = \Delta\lambda = (\Delta\mu + \Delta\nu) = \sigma + \rho.$$

Global Effectiveness Measure, *gem* = $-\tau$ = decrease in the infeasibility index.

## 4.1 Change in the primal infeasibility index ($\sigma = \Delta\mu$)

For each column j, the ratio $r_{ij} = \beta_i/\alpha_{ij}$ is calculated for all rows i = 1, 2, ..., m. It can be represented as in Figure-6.

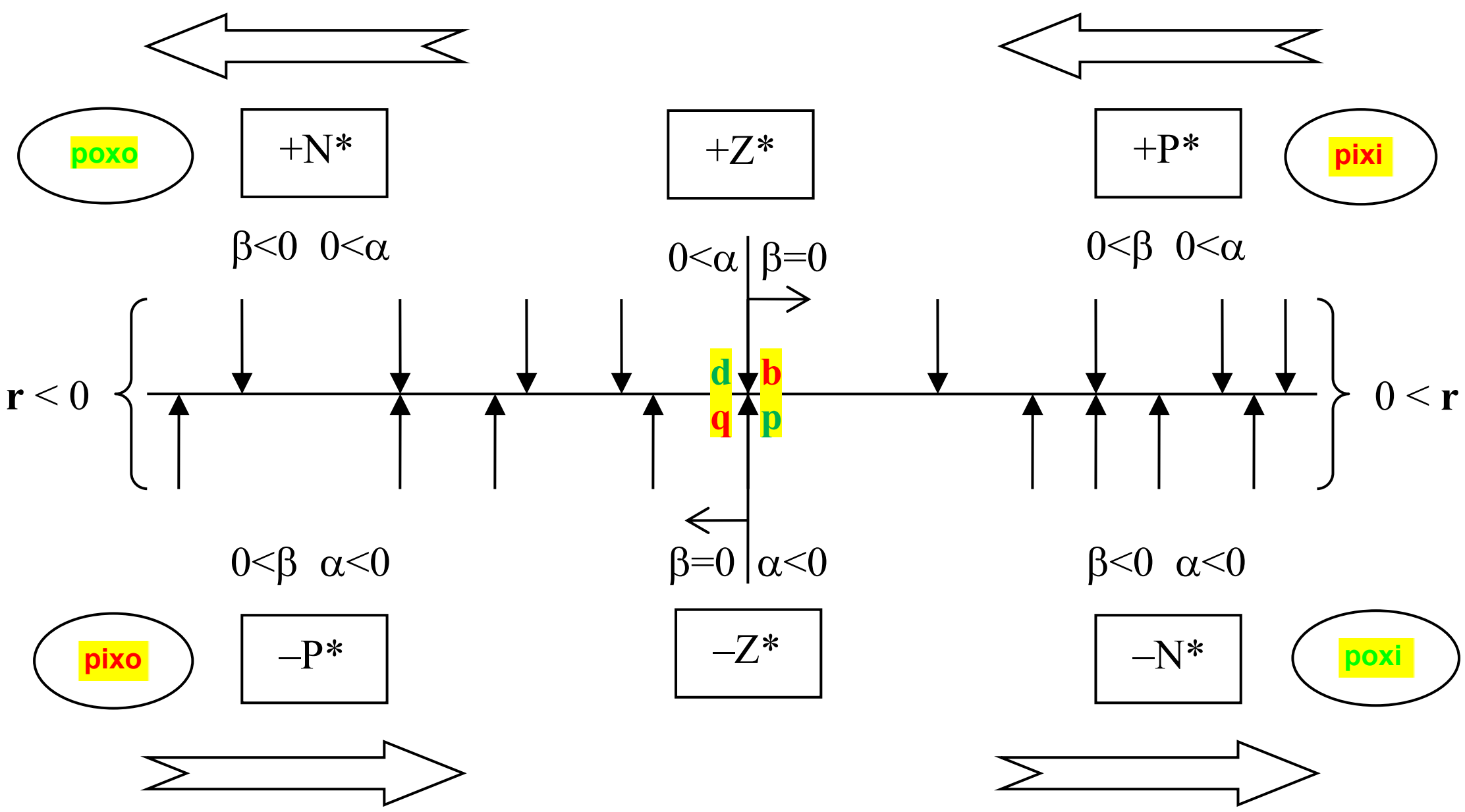

Figure-6: Ordered Pattern of $\{r_{iJ} = (\beta_i / \alpha_{iJ})\}$; i = 1, 2, ..., m

* => ('P' or 'N' or 'Z')

Let $\alpha_{IJ}$ be the chosen pivot element, then after pivoting, the values of $\beta$ are given by

$$r_{IJ} = \beta_I/\alpha_{IJ} \quad \text{and}$$

$$(\beta_I)' \leftarrow (\beta_I/\alpha_{IJ})$$

and $$(\beta_i)' \leftarrow \beta_i - (\beta_I/\alpha_{IJ})\alpha_{iJ}$$

The calculation of $\sigma$ can be divided into three different cases depending upon the value of $r_{IJ}$.

$\beta = 0$ is considered as feasible.

**(i)** $r_{IJ} = 0$ $(\beta_I = 0, \ \alpha_{IJ} \neq 0)$

Here, there will not be any change in the $\beta$ values, and therefore there will not be any change in the infeasibility index. Hence, $\sigma_{IJ} = 0$.

**(ii)** $0 < r_{IJ}$ $(0 \leq \beta_I, \ 0 < \alpha_{IJ} : \textit{pixi}$ -or- $\beta_I < 0, \ \alpha_{IJ} < 0 : \textit{poxi})$ : b–p

$r_{iJ} < 0$ will not affect the change in the infeasibility index. $r_{iJ} < 0$ can occur in two situations.

One of them is when $\beta_i < 0$ and $\alpha_{iJ} > 0$. In this case, it can be seen from the expression for $(\beta_i)'$ given above, that $(\beta_i)'$ will continue to be negative; hence the change in the infeasibility index is not affected.

The other is when $\beta_i > 0$ and $\alpha_{iJ} < 0$. In this case, it can be seen from the expression for $(\beta_i)'$ given above, that $(\beta_i)'$ will continue to be positive; hence the change in the infeasibility index is not affected.

$r_{iJ} > 0$ will affect the change in the infeasibility index. $r_{iJ} > 0$ can occur in two situations.

One is when $\beta_i \geq 0$ and $\alpha_{iJ} > 0$. In this case, it can be seen from the expression for $(\beta_i)$' given above, that $(\beta_i)$' will continue to be positive for the ratios $r_{iJ} > r_{IJ}$, $(\beta_i)$' will be zero for the ratios $r_{iJ} = r_{IJ}$, $(\beta_I)$' will continue to be positive, and $(\beta_i)$' will become negative for the ratios $r_{iJ} < r_{IJ}$. Hence the change in the infeasibility index is *increased* by the number of ratios $r_{iJ}$ which are less than $r_{IJ}$ and greater than or equal to zero.

The other is when $\beta_i < 0$ and $\alpha_{iJ} < 0$. In this case, it can be seen from the expression for $(\beta_i)$' given above, that $(\beta_i)$' will continue to be negative for the ratios $r_{iJ} > r_{IJ}$, $(\beta_i)$' will become zero for the ratios $r_{iJ} = r_{IJ}$, $(\beta_I)$' will become positive, and $(\beta_i)$' will become positive for the ratios $r_{iJ} < r_{IJ}$. Hence the change in the infeasibility index is *decreased* by the number of ratios $r_{iJ}$ which are less than or equal to $r_{IJ}$ and strictly greater than zero.

Therefore,

$b = \Sigma$ No. of rows i with $0 \leq \beta_i$, $0 < \alpha_{iJ}$ and $0 \leq r_{iJ} < r_{IJ}$ : {*pixi*}

$p = \Sigma$ No. of rows i with $\beta_i < 0$, $\alpha_{iJ} < 0$ and $0 < r_{iJ} \leq r_{IJ}$ : {*poxi*}

and

$\sigma_{IJ} = b - p$ : {*pixi-poxi*}

**(iii)** $r_{IJ} < 0$ ($0 \leq \beta_I$, $\alpha_{IJ} < 0$ : *pixo* -or- $\beta_I < 0$, $0 < \alpha_{IJ}$ : *poxo*) : 1+q–d

$r_{iJ} > 0$ will not affect the change in the infeasibility index. $r_{iJ} > 0$ can occur in two situations.

One of them is when $\beta_i \geq 0$ and $\alpha_{iJ} > 0$. In this case, it can be seen from the expression for $(\beta_i)$' given above, that $(\beta_i)$' will continue to be positive; hence the change in the infeasibility index is not affected.

The other is when $\beta_i < 0$ and $\alpha_{iJ} < 0$. In this case, it can be seen from the expression for $(\beta_i)$' given above, that $(\beta_i)$' will continue to be negative; hence the change in the infeasibility index is not affected.

$r_{iJ} < 0$ will affect the change in the infeasibility index. $r_{iJ} < 0$ can occur in two situations.

One is when $\beta_i \geq 0$ and $\alpha_{iJ} < 0$. In this case, it can be seen from the expression for $(\beta_i)$' given above, that $(\beta_i)$' will continue to be positive for the ratios $r_{iJ} < r_{IJ}$, $(\beta_i)$' will be zero for the ratios $r_{iJ} = r_{IJ}$, $(\beta_I)$' will become negative, and $(\beta_i)$' will become negative for the ratios $r_{iJ} > r_{IJ}$. Hence the change in the infeasibility index is increased by the number of ratios $r_{iJ}$ which are greater than $r_{IJ}$ and less than or equal to zero (one more, if $\beta_I$ is positive, since it will become negative).

The other is when $\beta_i < 0$ and $\alpha_{iJ} > 0$. In this case, it can be seen from the expression for $(\beta_i)$' given above, that $(\beta_i)$' will continue to be negative for the ratios $r_{iJ} < r_{IJ}$, $(\beta_i)$' will become zero for the ratios $r_{iJ} = r_{IJ}$, $(\beta_I)$' will become negative, and $(\beta_i)$' will become positive for the ratios $r_{iJ} > r_{IJ}$. Hence the change in the infeasibility index is decreased by the number of ratios $r_{iJ}$ which are greater than or equal to $r_{IJ}$ and strictly less than zero (one less, if $\beta_I$ is negative, since it will continue to be negative).

Therefore,

$q = \Sigma$ No. of rows i with $0 \leq \beta_i$, $\alpha_{iJ} < 0$ and $r_{IJ} < r_{iJ} \leq 0$ : {*pixo*}

$d = \Sigma$ No. of rows i with $\beta_i < 0$, $0 < \alpha_{iJ}$ and $r_{IJ} \leq r_{iJ} < 0$ : {*poxo*}

and

$\sigma_{IJ} = (q + 1) - d \quad$ if $(\beta_I > 0)$

or $\quad \sigma_{IJ} = q - (d - 1) \quad$ if $(\beta_I < 0)$

Hence we get,

$\sigma_{IJ} = 1 + q - d$ : {*1pixo-poxo*}

## 4.2 Change in the dual infeasibility index ($\rho = \Delta \nu$)

For each row i, the ratio $r_{ij} = \gamma_j/\alpha_{ij}$ is calculated for all columns j = 1, 2, ..., n. It can be represented as in Figure-7.

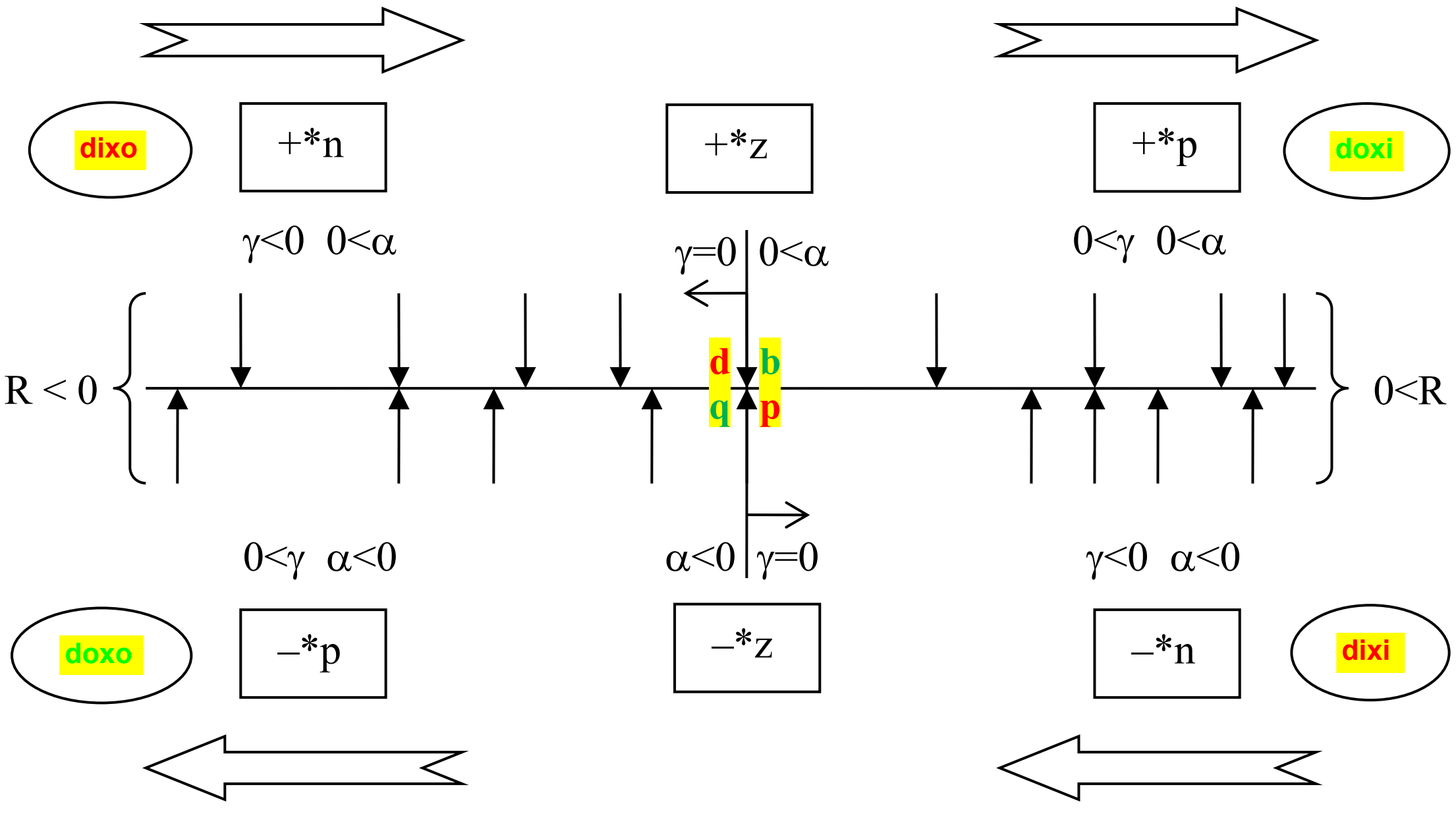

Figure-7: Ordered Pattern of $\{r_{Ij} = (\gamma_j / \alpha_{Ij})\}$; j = 1, 2, ..., n

* => ('p' or 'n' or 'z')

Let $\alpha_{IJ}$ be the chosen pivot element, then after pivoting, the values of $\gamma$ are given by

$$r_{IJ} = \gamma_J/\alpha_{IJ} \quad \text{and}$$

$$(\gamma_J)' \leftarrow -(\gamma_J/\alpha_{IJ})$$

and

$$(\gamma_j)' \leftarrow \gamma_j - (\gamma_J/\alpha_{IJ})\alpha_{Ij}$$

The calculation of $\rho$ can be divided into three different cases depending upon the value of $r_{IJ}$.

$\gamma = 0$ is considered as feasible.

**(i)** $r_{IJ} = 0 \qquad (\gamma_J = 0, \ \alpha_{IJ} \neq 0)$

Here there will not be any change in the $\gamma$ values, and therefore there will not be any change in the infeasibility index. Hence, $\rho_{IJ} = 0$.

**(ii)** $0 < r_{IJ}$ ($\gamma_J \leq 0$, $\alpha_{IJ} < 0$ : *dixi* -or- $\gamma_J > 0$, $\alpha_{IJ} > 0$ : *doxi*) : p–b

$r_{Ij} < 0$ will not affect the change in the infeasibility index. $r_{Ij} < 0$ can occur in two situations.

One of them is when $\gamma_j < 0$ and $\alpha_{Ij} > 0$. In this case, it can be seen from the expression for $(\gamma_j)'$ given above, that $(\gamma_j)'$ will continue to be negative; hence the change in the infeasibility index is not affected.

The other is when $\gamma_j > 0$ and $\alpha_{Ij} < 0$. In this case, it can be seen from the expression for $(\gamma_j)'$ given above, that $(\gamma_j)'$ will continue to be positive; hence the change in the infeasibility index is not affected.

$r_{Ij} > 0$ will affect the change in the infeasibility index. $r_{Ij} > 0$ can occur in two situations.

One is when $\gamma_j \leq 0$ and $\alpha_{Ij} < 0$. In this case, it can be seen from the expression for $(\gamma_j)'$ given above, that $(\gamma_j)'$ will continue to be negative for the ratios $r_{Ij} > r_{IJ}$, $(\gamma_j)'$ will be zero for the ratios $r_{Ij} = r_{IJ}$, $(\gamma_J)'$ will continue to be negative, and $(\gamma_j)'$ will become positive for the ratios $r_{Ij} < r_{IJ}$. Hence the change in the infeasibility index is *increased* by the number of ratios $r_{Ij}$ which are less than $r_{IJ}$ and greater than or equal to zero.

The other is when $\gamma_j > 0$ and $\alpha_{Ij} > 0$. In this case, it can be seen from the expression for $(\gamma_j)'$ given above, that $(\gamma_j)'$ will continue to be positive for the ratios $r_{Ij} > r_{IJ}$, $(\gamma_j)'$ will become zero for the ratios $r_{Ij} = r_{IJ}$, $(\gamma_J)'$ will become negative, and $(\gamma_j)'$ will become negative for the ratios $r_{Ij} < r_{IJ}$.
Hence the change in the infeasibility index is *decreased* by the number of ratios $r_{Ij}$ which are less than or equal to $r_{IJ}$ and strictly greater than zero.

Therefore,

$$p = \Sigma \text{ No. of columns j with } \gamma_j \leq 0,\ \alpha_{Ij} < 0 \text{ and } 0 \leq r_{Ij} < r_{IJ} \quad : \{dixi\}$$

$$b = \Sigma \text{ No. of columns j with } 0 < \gamma_j,\ 0 < \alpha_{Ij} \text{ and } 0 < r_{Ij} \leq r_{IJ} \quad : \{doxi\}$$

and

$$\rho_{IJ} = p - b \qquad : \{dixi\text{-}doxi\}$$

**(iii)** $r_{IJ} < 0$ ($\gamma_J \leq 0$, $\alpha_{IJ} > 0$ : *dixo* -or- $\gamma_J > 0$, $\alpha_{IJ} < 0$ : *doxo*) : 1+d–q

$r_{Ij} > 0$ will not affect the change in the infeasibility index. $r_{Ij} > 0$ can occur in two situations.

One of them is when $\gamma_j < 0$ and $\alpha_{Ij} < 0$. In this case, it can be seen from the expression for $(\gamma_j)$' given above, that $(\gamma_j)$' will continue to be negative; hence the change in the infeasibility index is not affected.

The other is when $\gamma_j > 0$ and $\alpha_{Ij} > 0$. In this case, it can be seen from the expression for $(\gamma_j)$' given above, that $(\gamma_j)$' will continue to be positive; hence the change in the infeasibility index is not affected.

$r_{Ij} <0$ will affect the change in the infeasibility index. $r_{Ij} < 0$ can occur in two situations.

One is when $\gamma_j \leq 0$ and $\alpha_{Ij} > 0$. In this case, it can be seen from the expression for $(\gamma_j)$' given above, that $(\gamma_j)$' will continue to be negative for the ratios $r_{Ij} < r_{IJ}$, $(\gamma_j)$' will be zero for the ratios $r_{Ij} = r_{IJ}$, $(\gamma_J)$' will become positive, and $(\gamma_j)$' will become positive for the ratios $r_{Ij} > r_{IJ}$. Hence the change in the infeasibility index is increased by the number of ratios $r_{Ij}$ which are greater than $r_{IJ}$ and less than or equal to zero (one more, if $\gamma_J$ is negative, since it will become positive).

The other is when $\gamma_j > 0$ and $\alpha_{Ij} < 0$. In this case, it can be seen from the expression for $(\gamma_j)$' given above, that $(\gamma_j)$' will continue to be positive for the ratios $r_{Ij} < r_{IJ}$, $(\gamma_j)$' will become zero for the ratios $r_{Ij} = r_{IJ}$, $(\gamma_J)$' will become positive, and $(\gamma_j)$' will become negative for the ratios $r_{Ij} > r_{IJ}$. Hence the change in the infeasibility index is decreased by the number of ratios $R_{Ij}$ which are greater than or equal to $r_{IJ}$ and strictly less than zero (one less, if $\gamma_J$ is positive, since it will continue to be positive).

Therefore,

$$d = \Sigma \text{ No. of columns j with } \gamma_j \leq 0,\ 0 < \alpha_{Ij} \text{ and } r_{IJ} < r_{Ij} \leq 0 : \{\textit{dixo}\}$$

$$q = \Sigma \text{ No. of columns j with } 0 < \gamma_j,\ \alpha_{Ij} < 0 \text{ and } r_{IJ} \leq r_{Ij} < 0 : \{\textit{doxo}\}$$

and

$$\rho_{IJ} = (d + 1) - q \quad \text{if} \quad (\gamma_J < 0)$$

or

$$\rho_{IJ} = d - (q - 1) \quad \text{if} \quad (\gamma_J > 0)$$

Hence we get,

$$\rho_{IJ} = 1 + d - q \qquad : \{\textit{1dixo-doxo}\}$$

## 5. Characterization of a Pivot Element/Cell

Each potential pivot element/cell in the Compact *Sym*metric Tableau is characterized by a "cell-type". The cell type of a cell in $I^{th}$ row and $J^{th}$ column consists of three components; they are – (1) the sign of $\alpha_{IJ}$ (either '0', '+' or '–'); (2) the sign of $\beta_I$ (either 'Z', 'P' or 'N'); and (3) the sign of $\gamma_J$ (either 'z', 'p' or 'n') – that is, [{0|+|-}{Z|P|N}{z|p|n}]. Hence, there are a total of 27 different possible cell types. However, if $\alpha$ is zero or numerically near-zero, it will not be a potential pivoting cell. Therefore, the nine cell types with $\alpha \approx 0$ may as well be combined together and the new cell type designated as [0**].

## 6. CST-Signature of a Compact *Sym*metric Tableau

We have designed a unique CST-*signature* that distinctly identifies any specific Goldman-Tucker Compact *Sym*metric Tableau associated with a given Linear Programming Problem. It is a string of length n+m (number of columns + number of rows). The first n entries are chosen from the character set {n, z, p} depending upon whether the $\gamma$ value is negative, zero or positive respectively; the next m entries are chosen from the character set {P, Z, N} depending upon whether the $\beta$ value is positive, zero or negative respectively. The positions of these entries in the signature string are fixed conveniently with respect to the initial tableau, as in considering the *lexicographic ordering* of these n+m parameters. As part of the pivot selection process, during the infeasibility analysis, corresponding to a specific possible choice of the pivot element, the CST-*signature* string can be generated for the anticipated resultant tableau. It can be compared with the signatures of all the previous tableaus stored in a dictionary, to facilitate the detection of any possible imminent cycle, since each CST-*tableau* has a unique CST-*signature* and each CST-*signature* uniquely identifies a CST-*tableau*. If & when an imminent cycle is detected, such pivot choice can be avoided and possible alternative pivot choice may be considered for actual pivoting operation. This methodology will be useful especially in situations of degeneracy. Note that an alternative binary string CST-signature may be defined by mapping the alphabet symbols {n,z,p} to 0 indicating the primal non-basic status, and mapping the alphabet symbols {P,Z,N} to 1 indicating primal basic status, of the corresponding parameters, as a framework to represent only the combinatorial information without any reference to the feasibility information.

## 7. Degeneracy

When no further decrease in the infeasibility index is indicated by every/all possible pivot selections (the best pivot choice corresponding to the largest possible decrease or non-increase in the infeasibility index is itself a case of no change in the infeasibility index) the situation may correspond to either of the two cases: (1) the tableau is indeed the terminal tableau, or (2) the best possible choice of the pivot element is in some degenerate row/column. In the first case, the *spdspds* algorithm terminates and the tableau may be analyzed to classify it into one of the six possible categories as described later. The presence of primal/dual degeneracy (either in the terminal tableau itself, as in the first case, or may even be some temporary/intermediate degeneracy encountered enroute towards the optimum/terminal tableau, as in the second case) may be identified by the presence of *zero* in some row/column. The situation wherein the best pivot choice is itself in some degenerate row/column, requires further analysis, before declaring that the tableau is indeed the terminal tableau - it may indeed be either a case of temporary/intermediate degeneracy or a case of degeneracy in a non-optimum but terminal tableau.

A pivot selection on any degenerate row/column with primal/dual degeneracy doesn't lead to any change in the primal/dual feasibility (refer to the algebra/arithmetic of the pivoting process presented earlier) although the coefficient matrix entries will surely be changed in the resulting tableau. Therefore, the infeasibility analysis needs to be carried out for the resulting tableau so as to check & confirm whether any further reduction in the infeasibility index is now attainable; if so, the usual iteration process continues; if not, then further pivoting on the degenerate rows/columns needs to be continued (as per a *predetermined lexicographically ordered sequence of the problem parameters*) among those possible pivot options while also keeping track of the CST-signature of the tableaus to detect & avoid possible impending cycling.

## 8. Algorithm Termination - Terminal Tableau Types

When further reduction in the infeasibility index is not possible, the tableau is checked & confirmed to be the *terminal tableau*. The possible terminal tableau types as shown in Figure-8 belong to one of the six possible categories listed here. Label F is for *feasible-basic-finite*, $\infty$ for *feasible basic-finite/non-basic-infinite* and $\Phi$ for *infeasible* - corresponding to both primal and dual *variables*.

This classification into *six categories* as presented here is a refinement over the well-known classical approach; one that enables to distinctly identify the case wherein primal (dual) has a *feasible basic-finite/non-basic-infinite* optimum with finite value for the objective function, while the dual (primal) has a *feasible-basic-finite* optimum. This distinction from the classical approach arises because we give primary emphasis on the classification based on the nature of the decision variables at termination, and give secondary emphasis on the finiteness (or otherwise) of the objective function value. The set of criteria (logical test for the existence and/or non-existence of specific cell types) for this classification scheme, to be used after ensuring that the tableau is indeed the *terminal tableau*, is also presented along with Figure-8.

| | D = F | D = ∞ | D = Φ |
|---|---|---|---|
| **P = F** | • Θ 0<br>⊕ * +<br>0 − • | + 0<br>* +<br>Θ • | |
| **P = ∞** | − * ⊕<br>0 − • | | • − ⊕<br>Θ + • |
| **P = Φ** | | • ⊕<br>+ −<br>Θ • | • • ⊕<br>• • −<br>Θ + • |

(1) {P=F,D=F}; (2) {P=F,D=∞}; (3) {P=∞,D=F}; (4) {P=∞,D=Φ}; (5) {P=Φ,D=∞}; (6) {P= Φ,D= Φ}

$$(P{=}F) : (\mu = 0)\ \&\ [\nexists\{(i\in R,\ j\in C)\} : \{(\alpha_{ij} < 0)(0 \le \beta_i)(0 \le \gamma_j)\}];$$

$$(P{=}\infty) : (\mu = 0)\ \&\ [\exists\{(i\in R,\ j\in C)\} : \{(\alpha_{ij} < 0)(0 \le \beta_i)(0 \le \gamma_j)\}];$$

$$(P{=}\Phi) : (\mu > 0)\ \&\ [\forall\{(I\in R,\ J\in C)\} : \{0 \le \sigma_{IJ}\}];$$

$$(D{=}F) : (\nu = 0)\ \&\ [\nexists\{(i\in R,\ j\in C)\} : \{(0 < \alpha_{ij})(\beta_i \le 0)(\gamma_j \le 0)\}];$$

$$(D{=}\infty) : (\nu = 0)\ \&\ [\exists\{(i\in R,\ j\in C)\} : \{(0 < \alpha_{ij})(\beta_i \le 0)(\gamma_j \le 0)\}];$$

$$(D{=}\Phi) : (\nu > 0)\ \&\ [\forall\{(I\in R,\ J\in C)\} : \{0 \le \rho_{IJ}\}];$$

Figure-8: Six Categories for Terminal Tableau Data Pattern.

− negative; 0 zero; + positive;

Θ non-positive; * any value; ⊕ non-negative; • un-analyzed

Note that the situations corresponding to primal/dual *degeneracy* with the resultant dual/primal *multiplicity* and also that of *infinite rays* are all easily discernible in this classification scheme.

In case of terminal infeasibility, a judicious use of either the *change in primal infeasibility index* ($\sigma=\Delta\mu$) or *change in dual infeasibility index* ($\rho=\Delta\nu$) instead of the *change in overall infeasibility index* ($\tau=\Delta\lambda$) can be utilized to arrive at an *almost primal feasible* or an *almost dual feasible* tableau, if one requires such an output for further problem analysis etc. This maneuverability can be adapted to the needs of the problem and that shows the *versatility* of the *spdspds* approach.

## 9. Symmetric Primal Dual *Sym*plex Pivot Selection Strategy

(1) The proposed *spdspds* approach can be used to solve any LP problem, by first converting it into the standard/canonical form before proceeding further. In performing such transformation, it is possible to enhance the overall efficiency by the following approach:

(a) Free/bounded variables can be replaced by non-negative variables not double in number, but *only just one extra in number.* Equations can be replaced by inequalities not double in number, but *only just one extra in number.*

(b) No need for use of artificial variables; the *initial basic solution* need not necessarily be feasible; and therefore no need to rely on two-phase method or big-M method, etc.

(2) Any nonzero element of the coefficient matrix in the tableau can be a potential pivot element. Note that the classical primal standard pivot (*doxipixi* pivot) with cell-type [+Pp] and the classical dual standard pivot (*poxidixi* pivot) with cell-type [–Nn] are the most often encountered ones, simply because of their classical appeal. Although for example, a *doxopixo* pivot with cell type [–Pp] or a *poxodixo* pivot with cell-type [+Nn] may usually seem to be simply unacceptable or extremely unlikely, it may indeed turn out that such a pivot choice can yield a surprisingly large decrease in the infeasibility index, depending on the tableau data entries, especially as an intermediary pivoting iteration, and therefore worth the consideration. The versatility of our spdspds algorithm allows for such pivot selection decisions.

(3) The proposed concept of *infeasibility index* along with the use of a *measure of goodness* for the pivot selection determined as/by the anticipated *decrease in the infeasibility index* considered as a *global effectiveness measure* (*gem* = -τ) arising due to the specific pivot selection, turns out to be a great grand breakthrough in achieving the ultimate performance challenge in the use of symplex method for solving linear programming problems. One can use either σ (if $m \leq n$) or ρ (if $n \leq m$) as a tie-breaking measure to choose among alternatives having the same value for τ. Also note that the *entire infeasibility analysis is independent of and unaffected by any scaling of rows/columns.*

(4) Problem of *cycling* enroute towards the terminal tableau, possibly caused by some *intermediate degeneracy* will get *prevented* because of the very nature of the *spdspds* pivot selection strategy, seeking the best possible decrease of the infeasibility index through a non-degenerate pivot. If and when encountered with a possible degeneracy that seems to block further progress in decreasing the infeasibility index, the *lexicographic ordering of the variables* may be used as a tie-breaking mechanism to guide the choice of pivot element, thus coursing through a part of the cycle, until a point is reached wherefrom *spdspds* finds a step down the infeasibility path again, unless it happens to be the terminal tableau. So, *spdspds* is effectively *immune* to problem of *cycling* caused by any possible *intermediate degeneracy.*

(5) Potential possibility of *numerical instability* can be avoided by careful elimination of poor choices of the pivot element, using *appropriate filters* - even if it requires going for a next-best pivot choice in terms of *gem* (*decrease in the infeasibility index*) - in order to avoid such treacherous pathways leading to numerically disastrous computational behavior - that is effectively like first exchanging I for K followed by an exchange of K for J in order to achieve an exchange of I for J - equivalent to using an appropriate *pre-conditioner* without actually doing so.

## 10. Computational Complexity

For a linear programming problem represented in its standard canonical form using the CST of size (m)x(n) the size(length) of the input data (problem size) $L$ can be taken as ~[(m+1)*(n+1)].

In the worst-case, the infeasibility index of the initial tableau will at most be {(m)+(n)} and at each *spdspds* symplex pivoting iteration the infeasibility index gets reduced by at least *one* so that it takes at most {(m)+(n)} *spdspds* symplex pivoting iterations to reach the optimum solution if one exists; or may even be well before that in order to report the infeasibility status of the given problem. Therefore, the number of *spdspds* symplex pivoting iterations required to solve the linear programming problem is, $O(L^{1/2})$ in the worst-case. However, the expected number of *spdspds* symplex pivoting iterations is ≤min{(m),(n)} because *every symplex pivot choice corresponds to a decrease in the infeasibility index and also that the basic/non-basic status of each variable is not expected to switch around much except in rare instances of potential numerical instability etc.*

Each *spdspds* symplex pivoting operation requires a complete *analysis of the infeasibility status* before performing the actual pivoting operation, and that itself requires ≤[2*(m)*(n)] floating point operations, that is $O(L)$. Each *spdspds* symplex pivoting operation requires ≤[3*(m+1)*(n+1)] floating point operations, which is again $O(L)$.

Thus, the total computational work involved in solving the linear programming problem is ≤{(m)+(n)}*[5*m*n + 3*m + 3*n] *floating-point arithmetic operations*, which is of course $O(L^{1.5})$ - that is *super-linear sub-quadratic polynomial time complexity bound.* We assert that this is indeed the lower limit in worst-case computational complexity for any symplex/simplex based algorithm for LP problem and that in fact it represents the computational complexity of LP problem itself.

It can be shown easily that the *space complexity* is only $O(L)$ - that is, linear complexity bound.

Note that *spdspds* has *finite termination* and provides a definite *output*, that is, either the optimum solution if & when one exists, or a point of *minimal infeasibility* beyond which no further *spdspds* symplex pivoting is possible towards any improvement in the feasibility of the given problem; which itself may provide insight as to the possible refinements in the problem formulation/modelling itself.

## 11. Directions for Future Research

0: The conceptual foundation for the design of *spdspds* may be expressed by the following reasoning:

$$\neg[(\exists I\in R)(\exists J\in C)]:[\{(\forall i\in R)(\forall j\in C)(i\neq I)(j\neq J)\{\tau_{IJ} < \tau_{ij} \leq 0\}].$$

The spdspds algorithm provides a lot of scope for future research work with a motivation to further enhance the computational efficiency - by *exploiting the complete information content* that is made available through the *infeasibility analysis* of initial tableau.

1: Analyze the entries in the initial tableau $T_0$ and make a statement regarding the basic/non-basic status of each of the variables in the optimum tableau $T_*$. This question is conceptually equivalent to asking, in the context of non-linear programming [16] as to whether a constraint will be *active* (non-basic) or otherwise at the optimum - although of course expecting an answer in that context may certainly not be practical. However, the situation with linear systems can be more promising so as to expect a possible attempt in answering that question - based only on a thorough analysis of either the initial tableau $T_0$, or equivalently, any tableau $T_k$ en-route towards the optimum - using a radically *redefined* concept of *binding/nonbinding* [17][18] *constraints* - implicitly achieved by the *spdspds* algorithm as explained herein below:

2: For any column $J \in C$ the set of rows $\{ineg\tau(J)\}$ with $\tau_{iJ} < 0$ that is $\{(i \in R):(\tau_{iJ}<0)]\}$ along with the associated *set of ratios* $\{r_{ineg\tau(J)}\}$ that is $\{\beta_{ineg\tau(J)}/\alpha_{ineg\tau(J).J}\}$ may be determined. Suppose there exists a row $I \in R$ having its ratio $\beta_I/\alpha_{IJ}$ for every column $J \in C$ that is farther away from zero beyond (relative to) the set of ratios $\{r_{ineg\tau(J)}\}$.

$$(\exists I \in R)(\forall J \in C)[\{\alpha_{IJ}=0\} \vee (\forall i \in R \backslash I)[(0 \leq \{\beta_{ineg\tau(J)}/\alpha_{ineg\tau(J).J}\} < \beta_I/\alpha_{IJ}) \vee (\beta_I/\alpha_{IJ} < \{\beta_{ineg\tau(J)}/\alpha_{ineg\tau(J).J}\} \leq 0)]]$$

Note that the feasibility/infeasibility of such a row (*non-binding?row*) will not be affected by any of our symplex pivoting operation.

Similarly, for any row $I \in R$ the set of columns $\{jneg\tau(I)\}$ with $\tau_{Ij} < 0$ that is, $\{(j \in C):(\tau_{Ij}<0)]\}$ along with the associated *set of ratios* $\{r_{jneg\tau(I)}\}$ that is $\{\gamma_{jneg\tau(I)}/\alpha_{I.jneg\tau(I)}\}$ may be determined. Suppose there exists a column $J \in C$ having its ratio $\gamma_J/\alpha_{IJ}$ for every row $I \in R$ that is farther away from zero beyond (relative to) the set of ratios $\{r_{jneg\tau(I)}\}$.

$$(\exists J \in C)(\forall I \in R)[\{\alpha_{IJ}=0\} \vee (\forall j \in C \backslash J)[(0 \leq \{\gamma_{jneg\tau(I)}/\alpha_{I.jneg\tau(I)}\} < \gamma_J/\alpha_{IJ}) \vee (\gamma_J/\alpha_{IJ} < \{\gamma_{jneg\tau(I)}/\alpha_{I.jneg\tau(I)}\} \leq 0)]]$$

Note that the feasibility/infeasibility of such a column (*non-binding?column*) will not be affected by any of our symplex pivoting operation.

This explains how *spdspds* algorithm effectively avoids any *unnecessary swapping back and forth of the basic/non-basic status of each variable* (that has been indeed the root cause of the severe computational inefficiency in the classical simplex method and all its variants - that has never been effectively addressed till now) by the very choice of the *spdspds* symplex pivot corresponding to the largest possible decrease in the infeasibility index at each symplex iteration.

Q: is there any further improvement possible beyond this? Suppose there exists a row i∈R having its ratio $\beta_i/\alpha_{ij}$ for every column j∈C that is closer to (no farther from) zero relative to the set of ratios $\{r_{ineg\tau(j)}\}$ - can we say that such a row (*binding*? *row*) will certainly be pushed to primal non-basic status by some pivoting operation? Suppose there exists a column j∈C having its ratio $\gamma_j/\alpha_{ij}$ for every row i∈R that is closer to (no farther from) zero relative to the set of ratios $\{r_{jneg\tau(i)}\}$ - can we say that such a column (*binding*? *column*) will certainly be pushed to primal basic status by some pivoting operation?

3: The combined effect of a simultaneous application of more than one symplex pivot (rather than a *sequence of symplex pivoting operations*) is indeed worth further detailed study, and can possibly lead to what may be called as the general "spdspds *omni*plex pivoting" operation.

4: The versatility of the spdspds algorithm allows for extensive combinatorial analysis along with numerical experimentation to explore further possibilities, including possible applications in *polyhedral combinatorics*, *network optimization*, etc. Also, for example, as in the classical approaches, one can possibly seek to first achieve primal or dual feasibility, more efficiently, by utilizing the best possible decrease in the corresponding component - either σ or ρ rather than τ, before going further with later iterations.

## 13. Conclusion

The proposed concept of *infeasibility index* is an *inverse measure of goodness* associated with a CST tableau. The *decrease in the infeasibility index* as a *global effectiveness measure* (*gem*) associated with each potential pivot element forms the basic foundation for the proposed *spdspds* algorithm. Also note that the *entire infeasibility analysis is independent of and unaffected by any scaling of rows/columns*. From an engineer's point of view, *scaling of the rows and columns* is equivalent to expressing the problem parameters (both dependent variables and independent variables) in a re-defined set of units-of-measurement, appropriately chosen, as a *preconditioning* step that may be required to improve the *numerical stability* of the overall computational process.

The proposed *spdspds* algorithm provides a novel viewpoint to the very same classic framework of Goldman-Tucker Compact-Symmetric-Tableau representation for LP; with a newly defined concept of *symplex pivot*; leading to an *efficient*, *robust* and *versatile* iterative solution strategy passing through a non-repetitive sequence of symplex tableaus requiring *minimum number of symplex iterations*, with a worst-case computational complexity of $O(L^{1.5})$.

The proposed concept of the *CST-signature* can be utilized to keep track of the computational path from the initial tableau to the terminal tableau, and for advance detection of possible cycling before it actually occurs; even though of course the *spdspds* algorithm is effectively immune to any cycling caused by possible intermediate degeneracy.

For the first time in the history of Linear Programming, we have achieved -

(0.0) effective utilization of the Goldman-Tucker Compact-*Sym*metric-Tableau (CST) which is a *unique symmetric* representation *common* to both the primal as well as the dual of a linear programming problem in its *standard canonical form*, and in which the algebra/arithmetic of the pivoting operation also gets represented by a unique single set of operations irrespective of whether it is a primal pivot or a dual pivot;

(1.0) with '*symplex*' a new avataar of 'simplex' in the 21st century, the concept of

(2.0) *infeasibility-index* for a P-D pair has been defined and effectively used as a

(3.0) *global-effectiveness-measure*, in the design of 'spdspds' algorithm, to guide the pivoting process, with *five guarantees* -

(3.1) improvement of the overall P-D feasibility at every iteration,

(3.2) immunity against cycling,

(3.3) finite termination,

(3.4) minimum number of spdspds symplex pivoting iterations,

(3.5) worst case computational complexity of $O(L^{1.5})$ -

(4.0) with a *re-defined concept of binding / non-binding constraints* -

(5.0) *implicitly incorporated* by this game-changer 'spdspds' algorithm, claiming to be the great grand breakthrough; having successfully resolved and reposed the *Linear Programming Performance Challenge of the millennium*; refer to the 9th among Smale's 18 "Mathematical Problems for the Next Century" [19];

(6.0) with consequent *immediate* as well as *lasting*, *deep* as well as *far-reaching* impact on the study, teaching and application of *Linear Programming* and *Algorithms* in particular as well as *Computational Complexity Analysis* in general.

## 14. Recommended Reading Materials

## 15. Acknowledgement

First and foremost, I wish to pay my humble respect and reverence to the pioneer George B. Dantzig, who introduced the (most versatile) simplex method (and its variants) for solving Linear Programming Problems, and can be rightly considered as the LP-Grand-Master forever.

I wish to put on record my gratitude towards Prof. Gerald L Thompson of CMU GSIA (now David A Tepper School of Business) who introduced the use of Goldman-Tucker Compact-*Sym*metric-Tableau for LP to me in his classes.

## 16. Declaration Regarding Affiliation and Funding

I, Dr(Prof) Keshava Prasad Halemane, hereby declare that I am a Professor retired as on 2017JAN31 from National Institute of Technology Karnataka Surathkal India, and I am not affiliated to any institution or organization or corporation or any other agency or whatever. This research work has been conducted entirely by me on my own as an Independent Researcher, and that I have not received any funding from any source other than my own savings, and I do not have any obligations or encumbrances of any kind, neither financial nor legal nor of any other kind, regarding the contents of this paper - of which I am the sole original author and creator. The work is based on algebraic/logical analysis; no data has been used; no conflict of interests.

## 17. Author Rights Retention Notice

## 18. DEDICATION

To my **ಅಜ್ಜ**(ajja) Karinja Halemane Keshava Bhat & **ಅಜ್ಜಿ**(ajji) Thirumaleshwari, to my **ಅಪ್ಪ**(appa) Shama Bhat & **ಅಮ್ಮ**(amma) Thirumaleshwari, for their *teachings through love*, that *quality matters more than quantity*; to my wife Vijayalakshmi for her *ever consistent love & support*; to my daughter Sriwidya.Bharati and my twin sons Sriwidya.Ramana & Sriwidya.Prawina for their *love & affection*.